\documentclass[12pt]{article}
\usepackage{amssymb}
\usepackage{amsfonts}
\usepackage{amsmath}
\newcommand{\R} {{\mathbb R}}
\newcommand{\Q} {{\mathbb Q}}
\newcommand{\Z} {{\mathbb Z}}

\newcommand{\N} {{\mathbb N}}
\newcommand{\C} {{\mathbb C}}

\begin{document}
\title{Problems on Bieberbach groups and flat manifolds}
\author{Andrzej Szczepa\'nski}
\date{\today}
\maketitle
\vskip 5mm
\hskip 60mm
{\em In memory of Charles B. Thomas}
\vskip 5mm
\begin{abstract}
We present about twenty conjectures, problems and questions about
flat manifolds. Many of them build the bridges between
the {\em flat world} and representation theory of the finite groups,
hyperbolic geometry and dynamical systems.
\end{abstract}
 
We shall present here some conjectures, problems and open questions related to 
flat manifolds. By flat manifold we understand a compact closed (without boundary)
Riemannian manifold with sectional curvature equal to zero. It is well known (see for example \cite{JW})
that any flat manifold $M^n$ can be considered as an orbit space $\R^n/\Gamma,$ where 
$\Gamma$ is a discrete, torsion free and cocompact subgroup of the group $E(n) = O(n)\ltimes \R^n =$ Isom$\R^n.$
It is easy to see that $\Gamma = \pi_{1}(M^n).$ If we remove the assumption that $\Gamma$ is torsion free 
we obtain a crystallographic group. Hence we also consider the problems
about them.
From the Bieberbach theorems (see \cite{LC}) any crystallographic group has a maximal free abelian subgroup (subgroup of all
translations) of finite index. Hence for any such group $\Gamma \subset E(n)$ we have a short
exact sequence
$$0\rightarrow \Z^n\rightarrow \Gamma \rightarrow H\rightarrow 0.$$
Here $H$ is a finite group which in case of torsion free crystallographic group (Bieberbach group)
is isomorphic to the holonomy group of manifold $\R^n/\Gamma.$  We call it the holonomy group of $\Gamma.$
Since the subgroup $\Z^n$ is maximal abelian, the {\it holonomy representation} $\phi_{\Gamma}:H\rightarrow GL(n,\Z),$
defined by conjugation in $\Gamma,$ is faithful. 
We shall denote by $\alpha$ an element of the second cohomology group $H^2(H,\Z^n)$ which is defined
by the above short exact sequence.
\vskip 1mm
\noindent
Finally we would like to mention that many conjectures, problems and open questions still make sense,
if the context of flat manifolds (or equivalently Bieberbach groups) is replaced 
by the context of almost flat manifolds (or equivalently almost Bieberbach groups).
Whenever this is the case, we shall denote this by $(\ast).$
\vskip 1mm
\noindent
We do not pretend completeness of our list. In fact, we do not touch
several questions about flat manifolds either
because of our ignorance about their status or because of insufficient interest from our point of view.
\vskip 1mm
\noindent
We thank K.Dekimpe and \fbox{Charles B. Thomas} for their advice and help.
We also would like to thank the referee for his (her) help and particulare for a correction of a few
errors and improving the language. 
\section{The classification problems}
It is well known (cf. \cite{LC}) that any finite group $G$ is a holonomy group of some flat manifold.
Hence we have.
\newtheorem{prob}{Problem}
\begin{prob}$^{\ast}$
Find the minimal dimension of a flat manifold with the holonomy group $G.$
\end{prob}
\newtheorem{rem}{Remark}
\begin{rem} 
The answer is known for cyclic groups (\cite{LC}, \cite{HH}), elementary abelian $p$ - groups,
dihedral groups, semidihedral groups, generalized quaternion groups (see \cite{HMSS}) and
simple groups $PSL(2,p)$ ($p$ is a prime number) \cite{WP}. There are also some other classes of finite
groups for which the minimal dimension is known. We do not pretend to give a complete list here.
Unfortunately, for the most finite groups
it seems to be very difficult question. For example, almost nothing is known for the symmetric groups.
The main obstruction is the calculation of the second cohomology of the finite group with special coefficients.
\end{rem}
\vskip 1mm
\begin{prob}
Do the calculation for those $p$-groups for which the cohomology is well-understood.
\end{prob}
\vskip 5mm
\noindent
An inductive classification of the Bieberbach groups follows from E. Calabi 
(\cite{EC}, \cite[Section 3.6 on page 124]{JW}).
This is a consequence of the fact that any subgroup of the torsion free crystallographic group
is again crystallographic. Hence any $n$-dimensional Bieberbach group $\Gamma$ with non trivial center
is an extension of some Bieberbach group $\Gamma'$ of dimension $n-1$ by integers $\Z$,
$$0\rightarrow \Gamma'\rightarrow \Gamma \rightarrow \Z \rightarrow 0.$$
Hence we have the following problem.
\begin{prob}$^{\ast}$
For a finite group $G$ find a Bieberbach group (with trivial center) of minimal
dimension.
\end{prob}
\begin{rem}  
The problem was considered in \cite{HS} (see also \cite{HMSS}, \cite{WP}),
where there are also some calculations. These include cyclic groups, 
elementary abelian $p$-groups (where $p$ is a prime number), dihedral $2$-groups,
semidihedral $2$-groups, the generalized quaternion $2$-groups and simple groups of Lie type. 
\end{rem}
\vskip 5mm
\noindent
In 1972 A. T. Vasquez \cite{AV} introduced, for any finite group $G,$ an invariant $n(G) \in \N.$ 
This is related to the class of flat manifolds with holonomy group $G.$
In 1989 in the beautiful paper \cite{CW} G. Cliff and A. Weiss , calculated $n(G)$ for any finite $p$-group. 
A purely algebraic definition of $n(G)$ is given in \cite{SII}. 
There is also given a characterization of the groups
with $n(G)=1.$
\begin{prob}
Calculate Vasquez invariant $n(G)$ for a finite group $G.$
\end{prob}
\begin{rem}  
There is a method, complementary to that of Calabi for the classification of Bieberbach groups with
given holonomy group $G.$ It says that any Bieberbach group $\Gamma$ with holonomy group $G$
and dimension $n \geq n(G)$ can be defined by a short exact sequence of groups
$$0\rightarrow \Z^{n-n(G)}\rightarrow \Gamma \rightarrow \Gamma_{G} \rightarrow 0,$$
where $\Gamma_{G}$ is a Bieberbach group of dimension $n(G)$, (cf. \cite{SII}).
\end{rem}
\vskip 5mm
\noindent
Let $f:M^n\rightarrow M^n$ be a continuous map, and
let $\tilde{f}:\R^n\rightarrow \R^n$ be its cover in the Euclidean space.
From the third Bieberbach theorem (see \cite{JW}), we have the following short exact sequence (cf. \cite{LC})
$$0\rightarrow Aff_0(M^n)\rightarrow Aff(M^n)\rightarrow Out(\Gamma)\rightarrow 0,$$
where 
$$Aff(M^n) = \{f:M^n \rightarrow M^n \mid \tilde{f}:\R^n\rightarrow \R^n \in GL(n,\R)\ltimes \R^n \}.$$
$Aff_0(M^n)$ denotes the connected component of the identity of the group of affine diffeomorphisms of
the flat manifold $M^n.$
It is isomorphic (see \cite{LC}) to 
$(S^1)^{\beta_1(M^n)},$
where $\beta_1$ denotes the first Betti number of the manifold, see 
\cite{HS}.
Moreover we have the following a short exact sequence (see \cite{HiS} and \cite{LC})
$$0\rightarrow H^1(H,\Z^n)\rightarrow Out(\Gamma) \rightarrow N_{\alpha}/H\rightarrow 0,$$
where $N = N_{GL(n,\Z)}(H)$ is the normalizer of the holonomy group $H$ in $GL(n,\Z)$ and
$$N_{\alpha} = \{n\in N\mid n\ast \alpha = \alpha \}.$$
Here $\ast$ is the standard action of the normalizer in cohomology, see \cite[page 168]{LC}.
Hence we have.
\newtheorem{theo}{Theorem}
\begin{theo} (\cite{SIV})
The following conditions are equivalent:
\vskip 1mm
(i) $Out(\Gamma)$ is a finite group,
\vskip 1mm
(ii) the normalizer $N$ is a finite group,
\vskip 1mm
(iii) the holonomy representation $\phi_{\Gamma}$
is $\Q$-multiplicity free and any 
\vskip 1mm
$\Q$-irreducible component is also $\R$-irreducible.
\end{theo}
{\bf Proof:} The equivalence of the conditions $(i)$ and $(ii)$ is obvious.
The proof of the last equivalence is much more difficult and we refer the reader to
\cite{SIV}. 
\vskip 5mm
\noindent
We can now fomulate the following,
\begin{prob} $^{\ast}$
Classify finite groups which are holonomy groups of flat manifolds $M^n = \R^n/\Gamma$
with finite outer automorphism group Out$\Gamma.$
\end{prob}
\begin{rem} Some progress was made in \cite{HiS}, e.g. for finite abelian groups.
\end{rem}
\vskip 1mm
\noindent
From a finite group
representation theory point of view, it might be interesting to consider
a similar classification for K\"ahler flat manifolds (see \cite{S}). Note that, here
$Out(\Gamma)$ is finite if and only if the holonomy representation is $\Q$-multiplicity free
and any $\R$-irreducible component is $\C$-reducible. 
\vskip 1mm
\noindent
The next conjecture is related to the last problem. 
\newtheorem{conj}{Conjecture} 
\begin{conj}
(\cite{SI})
For any finite group $H,$ there exist a Bieberbach group $\Gamma$ with $\Q$ multiplicity
free holonomy representation $\phi_{\Gamma}.$
\end{conj}
\begin{rem}  
If the first Betti number of the manifold $M^n$ is zero then,
from above, we have
$$Aff(M^n) \simeq Out(\Gamma).$$ This means that all information about the
symmetries of the manifild is built into the outer automorphism group. 
\end{rem}
\vskip 1mm
\noindent
It is interesting to compare the last problem with the work of M. Belolipetsky and A. Lubotzky 
\cite{BL} where they proved that for any finite group $G$ there exist a fundamental
group $\Gamma$ of some compact hyperbolic manifold such that $Out(\Gamma) = G.$
Hence we have.
\begin{prob}
Which finite groups $G$ occur as outer automorphism groups of Bieberbach groups with trivial center.
\end{prob}
In this connection we can ask if every finite subgroup of $GL(n,\Z)$ has a realisation as the centraliser of some
finite subgroup of $GL(n,\Z),$ for some $n.$ 
\begin{rem}  
In \cite{RW} there is an example of a flat manifold whose fundamental
group has trivial center and trivial outer automorphism group. 
\end{rem}
\section{The Anosov relation for flat manifolds and Entropy Conjecture}
Let $f:M\rightarrow M$ be a continuous map of a smooth manifold $M.$
In fixed point theory, there exist two numbers associated to $f$ which provide
some information on the fixed point set of $f:$ the Lefschetz number $L(f)$ and
the Nielsen number $N(f).$ It is known that the Nielsen number provides more
information about the fixed point set than $L(f)$ does, but $L(f)$ is easier
to calculate. For nilmanifolds however, D. Anosov \cite{An} showed that
$N(f) = \mid L(f)\mid$ for any self map $f$ of the nilmanifold. On the other hand,
he also observed that this result could not be extended to the class of all infranilmanifolds,
since he was already able to construct a counter-example on the Klein Bottle.
It was recently shown that for large families of flat manifolds (e.g. all flat
manifolds with an odd order holonomy group) the Anosov relation $N(f) = \mid L(f)\mid$ does
hold for any self map (see \cite{DRM} and \cite{DRP}).
It is therefore natural to ask.
\begin{prob}$^{\ast}$ 
(K. Dekimpe) Describe the class of flat manifolds (or more generally the class of infra-nilmanifolds)
on which the Anosov relation $$N(f) = \mid L(f)\mid$$ does hold for any map $f:M\rightarrow M.$
\end{prob}
One can also consider this problem from a different point of view and fix the flat
manifold $M$ and try to find all self maps $f$ for which the Anosov relation is valid.
\begin{prob}$^{\ast}$
(K. Dekimpe) Given a flat manifold $M.$ Is it possible to determine classes of
self maps $f:M\rightarrow M$ for which the Anosov relation $$N(f) = \mid L(f)\mid$$ does hold ?
\end{prob}
\vskip 3mm
\noindent
A first result in this direction can be found in \cite{MW}.
\vskip 5mm
\noindent
Let $n\in \Z^{+}$ and $\epsilon\in \R^{>0}$.  
Moreover let $(X,d)$ be a compact metric space and $f: X\rightarrow X$ a continuous map.
Put 
\begin{center}
$r(f,\epsilon)$ := lim$_{n\rightarrow \infty}$sup$\frac{1}{n}$ log max$\{\# Q: Q \subset X\},$
\end{center}
where $Q$ is such that
for any distinct points $x,y\in \Q$ 
\begin{center}
max$_{0\leq j\leq n-1}d(f^{j}(x),f^{j}(y)) \geq \epsilon.$
\end{center}
\vskip 1mm
\noindent
The {\em topological entropy} $h(f)$ is a nonnegative real number or $\infty$ defined
\begin{center}
$h(f)$ = lim$_{\epsilon \rightarrow 0}r(f,\epsilon)$ = sup$_{\epsilon \rightarrow 0}r(f,\epsilon).$
\end{center}
Now assume that $X$ is a compact smooth manifold $M$ of dimension $m.$
Let
$$H^{\ast}(f) : H^{\ast}(M,\R)\rightarrow H^{\ast}(M,\R)$$
be the linear map induced by $f$ on the cohomology space $H^{\ast}(M,\R):= \oplus_{i=0}^{m}H^{i}(M,\R)$
of $M$ with real coefficients.
By $sp(f)$ we denote the spectral radius of $H^{\ast}(f),$ which is a homotopy invariant.
In 1974 M.Shub asked in \cite[page 36]{MS} the extent to which the following inequality holds.
\begin{center}
(EC)\hskip 10mm log sp($f)\leq h(f).$
\end{center}
Then A. Katok in \cite{AK} put the following {\em Entropy Conjecture}.
\begin{conj}
(EC) holds for all continuous mappings $f:M\to M,$ if $M$ is a manifold with the universal cover
homeomorphic to $\R^n.$
\end{conj} 
It was proved for $M$ being tori and nilmanifolds, see \cite{MarP}.
\begin{prob}($^{\ast}$) Give a proof of the (EC) for flat manifolds.
\end{prob}
For any flat manifold $M$ with the first Betti number equal to zero
and finite outer automorphism group of the fundamental group $\pi_1(M)$ the above conjecture
is true. In fact, (see remark 5 and theorem 1) 
for any continuous map $f:M\to M$ the induced map $H^{\ast}(f)$ has a finite
order.
\vskip 1mm
\noindent
We send the reader to \cite{MarP} for more informations about the 
{\em Entropy Conjecture}$^{\ast}$.
\footnote{$\ast$}
\footnote{We would like to thank W. Marzantowicz for putting our attention on this problem.}
\section{Generalized Hantzsche-Wendt flat manifolds}
Let us agree to call any fundamental group of the $n$-dimensional 
flat manifold with holonomy group $(\Z_2)^{n-1}$ 
a generalized Hantzsche-Wendt $(GHW)$ Bieberbach group (see \cite{RS}).
If the manifold is orientable, we call the $GHW$ Bieberbach group orientable.
For dimension $3$ there is only one oriented flat manifold $M^3$ with holonomy group $\Z_2\oplus \Z_2.$
It was first considered by Hantzsche and Wendt. From the other side the fundamental group
$\pi_1(M^3)$ is group $F(2,6),$ 
where $F(r,n)$ is the group defined by the presentation
$$< a_{1}, a_{2}, ... , a_{n} \mid a_{1}a_{2}\cdots a_{r} = a_{r+1},
a_{2}a_{3}\cdots a_{r+1} = a_{r+2},  ... ,$$
$$a_{n-1}a_{n}a_{1}\cdots a_{r-2} = a_{r-1}, a_{n}a_{1}a_{2}\cdots a_{r-1} = a_{r} >,$$
where $r>0, n>0,$ and all subscripts are assumed to be reduced modulo $n$.
This group is named a Fibonacci group. 

\begin{prob}
Explain the relation between orientable GHW Bieberbach groups and the Fibonacci groups.
\end{prob}
\begin{rem} 
For any even natural number $2n$ there exists an
epimorphism 
$$\Phi_{2n} : F(2n,2(2n+1))\rightarrow \Gamma_{2n},$$
where $\Gamma_{2n}\subset E(2n+1)$ denote an orientable $GHW$ Bieberbach group
of rank $(2n+1),$ see \cite{SIII}. Nothing is known about 
the kernel of the epimorphism $\Phi_{2n},$ for $n > 1.$
\end{rem}
\begin{prob}
Give definition of any GHW Bieberbach group in terms of generators and relations.
\end{prob}
Let $\Gamma$ be a GHW Bieberbach group of dimension $n$ and with trivial center.
It is well known (\cite[Theorem 1.1, page 183]{FH}, \cite{JH}) 
that there is an epimorphism of $\Gamma$ onto the infinite
dihedral group $\Z_{2}\ast\Z_{2}.$
Hence we have a decomposition
$$(\ast)\hskip 5mm \Gamma \simeq \Gamma_{1}\ast_{X}\Gamma_{2},$$  
as the generalized free product of two Bieberbach groups of dimension $(n-1),$ amalgamated over
a subgroup $X$ of index two.
Moreover any GHW Bieberbach group ,with non-trivial center, is the semidirect product of $\Z$ with
some lower dimensional GHW Bieberbach group.
We have.
\begin{conj} For any GHW Bieberbach group $\Gamma$ with trivial center and dimension $n$, 
there exists a decomposition $(\ast)$ such that 
$\Gamma_1$ and $\Gamma_2$ are $GHW$ Bieberbach groups of dimension $n-1.$
\end{conj}
\begin{rem}  
The conjecture is true for $n\leq 4$ (see \cite[page 30 and 38]{JH}).
\end{rem}
For the remaining problem we have:
\vskip 1mm
\noindent
let $M^n$ be a flat oriented manifold with GHW fundamental Bieberbach group and dimension $n$.
It is well known, see \cite{RS}, that $M^n$ is a rational homology sphere. 
%and $M^n, n\geq 5$ does not have a spin structure.
\begin{prob}
What are the topological (geometrical) properties of $M^n$ ?
\end{prob}
\section{Flat manifolds and other geometries}
Let us recall the following well known question.
\vskip 1mm
\noindent
{\bf Question} 
(Farrell-Zdravkovska) \cite{FZ}) Let $M^n$ be a $n$-dimensional flat manifold. Is there a $(n+1)$-dimensional
hyperbolic manifold $W$ the boundary of which equals $M^n$ ?
\vskip 1mm
\noindent
We would like to mention that the above question is very close related to the following one.  
Let $V^{n+1}$ be a hyperbolic Riemannian manifold (constant negative curvature)
of finite volume. It is well known that $V^{n+1}$ has finite number of cusps and each cusp
is topologicaly $M^n\times \R^{\geq 0},$ where $M^n$ is $n$-dimensional flat manifold.
Is there some $V^{n+1}$ with one cusp
homeomorphic to $M^n\times \R^{\geq 0}$ ?
\vskip 1mm
\noindent
D.D.Long and A.Reid proved (\cite{LR}) that the answer to the problem above is negative if the $\eta$-invariant
of the signature operator of $M^{n}$ is not an integer.
To prove it they used the Atiyach, Patodi, Singer formula, where signature operator $D$ means some differential
elliptic operator and its $\eta$-invariant measure the symmetry of the spectrum of $D.$
They also proved that there
exists, already in dimension three, a flat manifold which signature operator has the $\eta$-invariant $\notin \Z.$
By our assumptions this method works only in dimension $4k-1.$ 
\begin{prob}
Classify flat manifolds with non-integral signature $\eta$-invariant.
\end{prob}
\vskip 1mm
\noindent
They also proved (\cite{LRI}) that any flat manifold has a "realization" as a cusp of
some hyperbolic orbifold of finite volume. 
\begin{rem}
We have the following flat manifolds of dimension two: the torus and the Klein bottle.
The first one has a "realization" as a cusp of the eight knot complement and the second one as a cusp of
the Gieseking manifold, (c.f. \cite{T}). 
In dimension three we have ten flat manifolds and each one has "realization" as a cusp of some
four dimensional hyperbolic manifold $W,$ which is not necessarly one-cusped, (see \cite{BN}).
\end{rem}
One major difficulty is the lack of examples of hyperbolic manifolds. In particular
we have.
\begin{prob}
Give an example of a four dimensional hyperbolic manifold of finite volume with only one cusp.
\end{prob}
\vskip 5mm
The $\eta$-invariant of the signature operator for a flat manifold is an important quantity as we
have observed above. Let us ask the same question about the
Dirac operator. We shall need our orientable flat manifold $M^n$ to have a spin structure.
It turns out that this is equivalent to the existence of a homomorphism $\epsilon:\Gamma \rightarrow Spin(n)$
such that the following diagram is commutative
$$\hskip 11mm Spin(n)$$
$$\epsilon \hskip 1mm \nearrow \hskip 10mm \downarrow \hskip 1mm \lambda_{n}$$
$$\Gamma \hskip 5mm \stackrel{r}{\rightarrow}\hskip 2mm SO(n),$$
where $\lambda_n$ is the universal covering and $r$ is the projection onto the linear part,
 \cite{DSS}.
\begin{prob} $^{\ast}$
Classify the holonomy groups of flat manifolds which admit a spin structure
\end{prob}
\begin{rem} 
It is known \cite[Proposition 1, Corollary 1]{DSS} that any flat manifold with holonomy group
of odd order or of order $2k$ , $k$ an odd number, has a spin structure.
However,
see \cite[Theorem 3.2]{MP}, for the holonomy group $\Z_2\times \Z_2,$ 
there exist orientable flat manifolds $M_1, M_2$ of dimension $6$ with different holonomy representations
$h_1, h_2: \Z_2\times \Z_2 \rightarrow GL(6,\Z),$ 
only one of which has a spin structure.
\end{rem}
\vskip 2mm
\noindent
{\bf Question}\hskip 2mm {\em
Is there an example of an oriented flat spin - manifold for which
the Dirac $\eta$-invariant is not equal to the signature $\eta$-invariant modulo $\Z$ ?}
\vskip 1mm
\noindent 
For the methods of calculation of the Dirac $\eta$-invariant see \cite{SS}.
\vskip 5mm
\noindent
A six dimensional flat manifold is Calabi-Yau if and only if
its holonomy representation  has the property that each $\R$-irreducible
summand, which is also $\C$-irreducible, occurs with even multiplicity.
Hence the classification of such manifolds is possible, see \cite{S}.
\begin{prob}
What are the properties of the Calabi-Yau flat manifolds. Do they have some
"mirror-symmetry" ?
\end{prob}
\vskip 5mm
\noindent
{\bf Acknowledgements}
\vskip 1mm
\noindent
This work was carried out in the framework of the project BIL01/C-31
for Bilateral Scientific Cooperation (Flanders--Poland) and was
supported by University of Gda\'nsk BW - 5100-5-0149-4.

Institute of Mathematics\\
University of Gda\'nsk\\
ul. Wita Stwosza 57\\
80-952 Gda\'nsk\\
Poland\\
E-mail: matas @ paula.univ.gda.pl
\end{document}